\documentclass [12pt]{article}
\usepackage{graphicx,amssymb,amsfonts,latexsym,amsmath,amsthm,times}
\usepackage{epsfig}
\usepackage{fancyhdr}
\usepackage{color}
\usepackage[english]{babel}
\usepackage[latin2]{inputenc}
\usepackage{graphics}
\usepackage{amsmath, amssymb, verbatim}
\usepackage{relsize}
\usepackage{graphicx}
\usepackage{subcaption}
\usepackage{epstopdf} 
\numberwithin{equation}{section}

\newtheorem{theorem}{Theorem}

\newtheorem{lemma}[theorem]{Lemma}

\newtheorem{proposition}[theorem]{Proposition}
\newtheorem{remark}[theorem]{Remark}

\begin{document}

%

\vspace*{2cm} \normalsize \centerline{\Large \bf Detailed analytic study of the compact pairwise model}
\normalsize \centerline{\Large \bf for SIS epidemic propagation on networks}
\vspace*{1cm}


\centerline{\bf Noémi Nagy, Péter L. Simon}

\vspace*{0.5cm}

\centerline{Institute of Mathematics, E\"otv\"os Loránd University Budapest, and}  \centerline{Numerical Analysis and Large Networks Research Group, Hungarian Academy of Sciences, Hungary}


%
%

\begin{abstract}
The global behaviour of the compact pairwise approximation of SIS epidemic propagation on networks is studied. It is shown that the system can be reduced to two equations enabling us to carry out a detailed study of the dynamic properties of the solutions. It is proved that transcritical bifurcation occurs in the system at $\tau = \tau _c =  \frac{\gamma n}{\langle n^{2}\rangle-n}$, where $\tau$ and $\gamma$ are infection and recovery rates, respectively, $n$ is the average degree of the network and $\langle n^{2}\rangle$ is the second moment of the degree distribution. For subcritical values of $\tau$ the disease-free steady state is stable, while for supercritical values a unique stable endemic equilibrium appears. We also prove that for subcritical values of $\tau$ the disease-free steady state is globally stable under certain assumptions on the graph that cover a wide class of networks.
\end{abstract}




\vspace*{1cm} \setcounter{equation}{0}
\section{Introduction}

Spreading processes on networks are widely studied by using stochastic and dynamical approaches \cite{MillerKissSimon, pastor2001PRE, porter}.
The mathematical model of such a process, like epidemic propagation on a graph, can
be formulated as a large system of linear ordinary differential equations.
The mathematical model is given by a graph the nodes of which can be in
different states, in the case of epidemic dynamics each node can be either susceptible
or infected. The state of the network containing $N$ nodes is given by an $N$-tuple of $S$ and
$I$ symbols, i.e. there are $2^N$ states altogether. The transition rules determine how the state
 of the network evolves by infection from one node to its neighbours and by recovery, when an $I$
  node becomes $S$ again. The system of master equations is formulated in terms of the probabilities
  of the states, i.e. the system consists of $2^N$ differential equations. Similar differential equations
   can be derived for modeling neural activity in a neural network (when each neurone, a node of the
    network, can be either active or inactive) or for the voter model describing the collective behaviour of voters.

Despite of the fact, that the mathematical model is relatively simple,
the analytical or numerical study of the system can be carried out only
 for small graphs or graphs with many symmetries like the complete graph
  or the star graph. For real-world large graphs the master equations are beyond tractability,
   hence the system is approximated by simple non-linear differential equations, called mean-field equations.
   These differential equations are formulated in terms of population level quantities, such as the expected
    number of susceptible and infected nodes, denoted by $[S]$ and $[I]$, or the average number of edges connecting
     different types of nodes, e.g. $[SI]$. The simplest mean-field model is written at the node level and the
     closure is applied to pairs. For SIS epidemic propagation the exact, population level differential equation
     takes the form
$$
\dot{[I]} =\tau [SI] -\gamma [I],
$$
where $\tau$ and $\gamma$ are nonnegative parameters, called infection and recovery rates. This equation becomes self-contained with a closure relation expressing $[SI]$ in terms of $[I]$. The simplest closure, expressing that susceptible and infected nodes are distributed randomly in the network, takes the form $[SI]\approx n[S][I]/N$, where $n$ is the average degree of the network and $[S]=N-[I]$. This approximation yields reasonable accuracy only under certain conditions, hence more accurate mean-field models, called pairwise models are formulated in terms of singles (nodes) and pairs (edges) and the closure is applied at the level of triples. The simplest one of these models is the homogeneous pairwise model:
\begin{eqnarray*}
\dot{[I]}&=&\tau [SI] -\gamma [I],  \\
\dot{[SI]}&=&\gamma ([II]-[SI])+\tau ([SSI]-[ISI]-[SI]), \\
\dot{[SS]}&=&2\gamma [SI]-2 \tau [SSI] , \\
\dot{[II]}&=&2\tau [SI]-2\gamma [II]+2\tau [ISI] ,
\end{eqnarray*}
in which  the triple closure $[ABC]=\frac{n-1}{n} [AB][BC]/[B]$ is used. These models, written in terms of singles and pairs, are widely applied since the early work of Matsuda et al. \cite{matsuda1992statistical} and Keeling et al. \cite{keeling1997}, and their unclosed forms were derived from exact master equations in \cite{lumpingcikk, taylor2012JMB}. This type of coarse-graining is not satisfactory when the network is strongly heterogeneous, i.e. there are nodes with low and high degree. Then, instead of using $[I](t)$ the average number of infected nodes in the network, the differential equations are formulated in terms of new variables, such as $[S_k](t)$ and $[I_k](t)$ which denote the average number of susceptible and infected nodes of degree $k$. Using these variables the closure relations can be defined in a more accurate way, hence these models perform significantly better in the case of heterogeneous networks. These models, called heterogeneous mean-field (or degree-based mean-field) and heterogeneous pairwise models, were introduced \cite{eames2002, pastor2001PRE}.

The heterogeneous pairwise model yields excellent approximation for Configuration Model random graphs, however, its size is still too large for analytical investigations. Hence, a reduced model, the compact pairwise model was developed in \cite{house2011insights}. In order to formulate the differential equations of the model let us denote by $N_{l}$ the number of nodes of degree ${n_{l}}$, $l=1,2,\ldots , L$ for a graph with $N$ nodes. The notations $ \langle n \rangle =\frac{1}{N}\sum_{l=1}^{L}n_l N_l$ and $\langle n^{2} \rangle =\frac{1}{N}\sum_{l=1}^{L}n_l^{2} N_l$ are used for the average degree and for the second moment of the degree distribution. The most important quantities under investigation are the average number of nodes in a given state with a given degree at time $t$ that are denoted by $[S_l](t)$ and $[I_l](t)$ for susceptible and infected nodes of degree $n_l$. For the average number of
$SI$, $SS$ and $II$ edges at time $t$, the notations $[SI](t)$, $[SS](t)$ and $[II](t)$ are applied, respectively. Then, the compact pairwise model takes the form
\begin{eqnarray}
 \dot{[S_l]}&=&\gamma [I_l]-\tau n_l [S_l] \frac{[SI]}{S_s}, \quad l=1,\ldots, L \label{alap_CP1}\\
 \dot{[I_l]}&=&\tau n_l [S_l] \frac{[SI]}{S_s} -\gamma [I_l],  \quad l=1,\ldots,L
 \label{alap_CP2}\\
 \dot{[SI]}&=&\gamma ([II]-[SI])+\tau ([SS]-[SI])[SI] Q_{CP}-\tau [SI], \label{alap_CP3}\\
 \dot{[SS]}&=&2\gamma [SI]-2 \tau [SS] [SI] Q_{CP}, \label{alap_CP4}\\
 \dot{[II]}&=&2\tau [SI]-2\gamma [II]+2\tau [SI]^2 Q_{CP} \label{alap_CP5},
\end{eqnarray}
where \[S_s=\sum_{l=1}^L n_l[S_l]=[SS]+[SI], \qquad Q_{CP}=\frac{1}{S_s^{2}} \sum_{l=1}^L (n_l-1)n_l [S_l].\]

While the dynamical behaviour and bifurcation analysis of the homogeneous mean-field and pairwise models
are considered to be folklore (nevertheless, their detailed study can be found in the textbook \cite{MillerKissSimon}),
 the detailed study of steady states and the global behaviour has not been carried out yet for the compact pairwise model.
 The aim of this paper is to investigate the dynamic of system \eqref{alap_CP1}-\eqref{alap_CP5} in detail. We will prove
  that there is a disease-free steady state (without infection) of the system and characterize its stability. It loses
   stability via transcritical bifurcation giving rise to an endemic steady state as the infection rate $\tau$ surpasses
    the critical value $\tau_c=\gamma \frac{ \langle n \rangle }{ \langle n^2 \rangle - \langle n \rangle }$. It will be
     shown that the endemic steady state is unique and applying a general transcritical bifurcation theorem, it is shown
      that the unique endemic steady state is asymptotically stable. Finally,
the global stability of the disease-free steady state is proved for a wide range of parameter
values.

\vspace*{0.5cm} \setcounter{equation}{0}
\section{Number of steady states}

In this section the equilibrium points of the compact pairwise model are studied. Adding the differential equations in system \eqref{alap_CP1}-\eqref{alap_CP5} it is obvious that the following proposition holds.
\begin{proposition} \label{prop1}
The conservation of singles $[S_l]+[I_l]=N_l$, for $l=1,\ldots, L$ and the conservation of  pairs $[SS]+2[SI]+[II]=nN$ ($S_s=[SS]+[SI]$) hold in system (\ref{alap_CP1})-(\ref{alap_CP5}), where $n=\langle n \rangle$ is the average degree.
\end{proposition}
As a consequence, the size of the system (\ref{alap_CP1})-(\ref{alap_CP5}) can be reduced in different ways by expressing some variables in terms of the others, as it will be shown later.

It is easy to see that system (\ref{alap_CP1})-(\ref{alap_CP5}) has a disease free steady state for any value of the parameter $\tau$, namely $[S_l]=N_l$, $[I_l]=0$, $l=1,\ldots,L$, $[SI]=0$, $[SS]=nN$, $[II]=0$. It will be verified that there is a critical value $\tau_{c}$, at which the system behaviour changes and another equilibrium point appears, which is called the endemic steady state. To prove this, equation \eqref{alap_CP2} will be omitted from the system. Then, by using the conservation of singles $[S_l]+[I_l]=N_l$, equation \eqref{alap_CP1} takes the form
\begin{eqnarray}
 \dot{[S_l]}&=&\gamma (N_l-[S_l])-\tau n_l [S_l] \frac{[SI]}{S_s}, \quad l=1,\ldots, L . \label{CP111}
\end{eqnarray}
Concerning the equilibria we have the following result.

\begin{theorem} \label{theorem1}
The critical value of the compact pairwise model is
$$
\tau_c=\gamma \frac{ \langle n \rangle }{ \langle n^2 \rangle - \langle n \rangle }.
$$
If $\tau_c<\tau$, then the model (\ref{CP111})-(\ref{alap_CP3})-(\ref{alap_CP5}) has a unique endemic steady state (with $[S_l]<N_l$).
\end{theorem}

\noindent {\sc Proof.} Let us denote the endemic steady state values of the variables $[S_l]$, $[SI]$, $[SS]$ and $[II]$ by $X_l$, $Z$, $U$, $V$, respectively. Putting zero in the left hand side of (\ref{CP111}) yields the following equation for $X_l$.

\begin{equation}
\gamma N_l=X_l \Big ( \gamma+ \tau n_l \frac{Z}{Z+U} \Big), \quad l=1, \ldots, L. \label{Xl}
\end{equation}
Expressing $X_l$, multiplying the equations by $n_l(n_l-1)$, and summing them for $l=1,\ldots, L$, we obtain
\begin{equation}\label{equ1}
\gamma \sum_{l=1}^{L}\frac{n_l(n_l-1)N_l}{\gamma +\tau n_l
\frac{Z}{Z+U}}=\sum_{l=1}^{L}n_l(n_l-1)X_l=Q_{CP}(U+Z)^2 \ .
\end{equation}
According to equation (\ref{alap_CP4}), $\tau Q_{CP}U= \gamma$ holds at the equilibrium point. Hence dividing equation (\ref{equ1}) by $Q_{CP}(U+Z)^2$, we get
$$
1=\frac{\tau U}{U+Z}\sum_{l=1}^{L}\frac{n_l(n_l-1)N_l}{\gamma (U+Z)+\tau n_l Z}.
$$
Equation (\ref{alap_CP5}) yields that in the steady state we have $\gamma V -\tau Z=\tau Z^2 Q_{CP}$, furthermore recalling $\tau Q_{CP}U= \gamma$ and $[SS]_c+2[SI]_c+[II]_c=nN$ leads to the following relation between $U$ and $Z$
\begin{equation}
\begin{split} \label{U_Z}
\gamma nNU=\gamma Z^2 +ZU(\tau+2\gamma)+\gamma U^2,
\end{split}
\end{equation}
which can be solved for $Z$ in terms of $U$. It is easy to see that for any $U\in [0,nN]$, equation (\ref{U_Z}) has a unique nonnegative solution for $Z$. Let this solution be denoted by $Z=g(U)$. We note that $g(0)=0$ and $g(nN)=0$.

A positive auxiliary function $f$ is defined as follows.
\[f(U):=\frac{\tau U}{U+g(U)}\sum_{l=1}^{L}\frac{n_l(n_l-1)N_l}{\gamma (U+g(U))+\tau n_l g(U)}=\frac{\tau U}{(U+g(U))^{2}} \sum_{l=1}^{L} \frac{n_l(n_l-1)N_l}{\gamma +\tau n_l \frac{g(U)}{(U+g(U))}}.\]
Then the existence and the uniqueness of the endemic steady state is equivalent to the fact that there is a unique  $U\in (0,nN)$ satisfying $f(U)=1$, since $U$ determines $Z$ via $Z=g(U)$, $X_l$ via \eqref{Xl} and $V$ via the conservation of pairs given in Proposition \ref{prop1}.

We prove the uniqueness of $U$ in three steps. First we verify that the limit of $f$ as $U\rightarrow 0$ is less than $1$. (Note that $f(0)$ is not defined.) To see this, we rearrange (\ref{U_Z}) as
\begin{equation}
\begin{split} \label{U_Z2}
\gamma nNU=\gamma (Z+U)^2+\tau ZU,
\end{split}
\end{equation}
leading to
\[\frac{U}{(U+Z)^2}=\frac{\gamma}{\gamma nN-\tau Z} . \]
Taking the limit $U\rightarrow 0$ yields
\[\lim _{U\rightarrow 0}\frac{U}{(U+g(U))^2}=\frac{1}{nN}.\]
Dividing equation  (\ref{U_Z}) by $Z$, one obtains
\[\lim _{U\rightarrow 0}\frac{U}{g(U)}=0, \textrm{ hence } \lim _{U\rightarrow 0}\frac{g(U)}{U+g(U)}=1. \]
Using these limits we obtain
\[ \lim_{U\rightarrow 0} f(U)=\frac{\tau}{nN}\sum_{l=1}^{L}\frac{n_l(n_l-1)N_l}{\gamma+\tau n_l}<\frac{\tau}{nN}\sum_{l=1}^{L}\frac{n_l(n_l-1)N_l}{\tau n_l}=\frac{N(n-1)}{Nn} <1. \]

Secondly, we show that $f(nN)>1$.
\[f(nN)=\frac{\tau}{nN} \sum_{l=1}^{L}\frac{n_l(n_l-1)N_l}{\gamma}=\frac{\tau}{\gamma} \frac{ \langle n^2 \rangle - \langle n \rangle }{ \langle n \rangle } =\frac{\tau}{\tau_c} . \]
Then, it is obvious, that $f(nN)>1$, when $\tau_c<\tau$.

The final step is to justify that the function $f$ is strictly increasing and continuous in the interval $(0,nN)$. The continuity of $f$ follows directly from the fact that the function $g$ is continuous in the interval $(0,nN)$. Hence it remained to show that $f$ is strictly increasing. To prove this, it is enough to guarantee that each term of the sum in the definition of $f$ is strictly increasing, namely the functions
\[U\mapsto\frac{\tau U}{U+g(U)}\frac{n_l(n_l-1)N_l}{\gamma (U+g(U))+\tau n_l g(U)}\]
need to be proved to be strictly increasing for all $l=1,\ldots,L$. The monotonicity of these functions is equivalent to the monotonicity of the functions
\[U\mapsto\frac{U}{U+g(U)}\frac{1}{\gamma (U+g(U))+\tau n_l g(U)}.\]
Using $\gamma (U+g(U))^2=\gamma nNU-\tau U g(U)$ again (see (\ref{U_Z2})) we get

\[\frac{U}{U+g(U)}\frac{1}{\gamma (U+g(U))+\tau n_l g(U)}=\frac{U}{(\gamma +\tau n_l)(U+g(U))^{2}-\tau n_l U (U+g(U))}=\]
\[=\frac{U}{\frac{\gamma +\tau n_l}{\gamma}(\gamma nNU-\tau U g(U))-\tau n_l U (U+g(U))}= \]
\[=\frac{1}{(\gamma+\tau n_l)nN-\tau n_l U-g(U)\tau(1+\frac{\tau}{\gamma}n_l+n_l)}.\]
To conclude, we prove that the function
\[U\mapsto (\gamma+\tau n_l)nN-\tau n_l U-g(U)\tau(1+\frac{\tau}{\gamma}n_l+n_l)\]
is a decreasing function in the interval $(0,nN)$. Omitting the positive constants, it is sufficient to show that the function
\[h_l(U):= U+g(U)(\frac{1}{n_l}+\frac{\tau}{\gamma}+1)\]
is strictly increasing in the interval $(0,nN)$. To this end, we will prove below that $h_l'(U)=1+g'(U)(\frac{1}{n_l}+\frac{\tau}{\gamma}+1)>0$ when $U\in (0,nN)$.

Consider first the  explicit formula of the function $g$, given by the solution of the quadratic equation (\ref{U_Z}) as
\[g(U)=\frac{1}{2\gamma}\big(-U(\tau+2\gamma)+\sqrt{U^2(\tau^2+4\gamma \tau)+4UnN \gamma^2} ~\big) .\]
Then the first and the second derivatives of $g$ take the form
\[g'(U)=\frac{1}{2\gamma}\Big(-(\tau+2\gamma)+\frac{1}{2}\frac{2U(\tau^2+4\gamma \tau)+4nN \gamma^2}{\sqrt{U^2(\tau^2+4\gamma \tau)+4UnN \gamma^2}} ~\Big),\]
\[g''(U)=\frac{-2n^2N^2 \gamma^3}{(U^2(\tau^2+4\gamma \tau)+4UnN \gamma^2)^{\frac{3}{2}}} <0.\]

It is obvious, that $g''$ is negative in the interval $(0,nN)$, thus the derivative $g'$ is strictly decreasing. Note that $g'(nN)=-\frac{\gamma}{\tau+2\gamma}<0$ and $\lim_{U\rightarrow 0} g'(U)=+\infty$, hence the derivative $g'$ has a unique root in the interval $(0,nN)$. In the subinterval, where $g'\geq 0$, the function $h_l'$ is positive, thus only that subinterval should be investigated, where $g'<0$.

Because of the monotonicity of $g'$ the following inequalities hold in the interval, where $g'<0$:
$$
h_l'(U)=1+g'(U)(\frac{1}{n_l}+\frac{\tau}{\gamma}+1)>1+g'(U)(1+\frac{\tau}{\gamma}+1)>1-\frac{\gamma}{\tau+2\gamma}(2+\frac{\tau}{\gamma})=0.
$$
Thus $h_l'(U)$ is positive for every $U\in(0,nN)$, which completes the proof.  $\Box$

\vspace*{0.5cm} \setcounter{equation}{0}
\section{Stability of the steady states } \label{seceSI}

The stability of the disease free equilibrium is studied first. In order to carry out the stability analysis, the size of system (\ref{alap_CP1})-(\ref{alap_CP5}) is reduced. According to Proposition \ref{prop1}, one of the equations for the singles and one of those for the pairs can be omitted. Omitting equations (\ref{alap_CP2}) and (\ref{alap_CP4}) leads to the following reduced form of the compact pairwise system.
\begin{eqnarray}
 \dot{[S_l]}&=&\gamma (N_l-[S_l])-\tau n_l [S_l] \frac{[SI]}{S_s}, \quad l=1,\ldots, L \label{CPr12ppp}\\
 \dot{[SI]}&=&\gamma ([II]-[SI])+\tau (nN-3[SI]-[II])[SI] Q_{CP}-\tau [SI], \label{CPr3ppp}\\
\dot{[II]}&=&2\tau [SI]-2\gamma [II]+2\tau [SI]^2 ,
Q_{CP},\label{CPr4pppp}
\end{eqnarray}
with
\[S_s=\sum_{l=1}^L n_l[S_l], \qquad Q_{CP}=\frac{1}{S_s^{2}} \sum_{l=1}^L (n_l-1)n_l [S_l].\]

\begin{proposition} \label{prop2}
In the compact pairwise model the stability of the disease free steady state changes at the critical value $\tau_{c}$. For $\tau<\tau_{c}$ the disease free steady state is asymptotically stable and for $\tau_{c}<\tau$ it is unstable.
\end{proposition}
\noindent {\sc Proof.}

Let us consider the linearisation of the right hand side of system (\ref{CPr12ppp})-(\ref{CPr4pppp}). The Jacobian matrix $J$ at the disease free equilibrium takes the form
$$
J=\begin{pmatrix}-\gamma I&R\\
                      0 & P
\end{pmatrix},
$$
where $0$ is a $2$-by-$L$ zero matrix and $I$ is a $L$-by-$L$ identity matrix, $R$ denotes a $L$-by-$2$ matrix
and $P$ is the $2$-by-$2$ matrix
$$
P=\begin{pmatrix}\alpha &\gamma\\
                      2\tau & -2\gamma
\end{pmatrix}
$$
with $ \alpha=\tau \frac{ \langle n^2 \rangle - \langle n \rangle }{ \langle n \rangle }-(\tau+\gamma)$.

The block structure of matrix $J$ implies that the multiplicity of the eigenvalue $-\gamma $ is $L$ and the other two eigenvalues of $J$ are the roots of the polynomial
\begin{equation}\label{poli}
\lambda^{2}+\lambda(2\gamma-\alpha)-2\gamma(\alpha+\tau).
\end{equation}

The disease free steady state is asymptotically stable, if all the eigenvalues have negative real part, which is satisfied when the coefficients of the quadratic polynomial in \eqref{poli} are positive, that is the relations $\alpha+\tau<0$ and $2\gamma-\alpha>0$ hold. The first relation implies the second one, hence the disease free steady state is asymptotically stable when $\alpha+\tau<0$, which is equivalent to the inequality $\tau<\tau_{c}$. Otherwise, if $\tau_{c}<\tau$, a positive eigenvalue appears and the disease free steady state becomes unstable. $\Box$

We have seen in Theorem \ref{theorem1} that the endemic equilibrium appears at the same critical value of $\tau$ where the disease-free steady state loses its stability. This observation suggests that transcritical bifurcation occurs at the critical value $\tau_{c}$. To investigate the stability of the endemic steady state, the following general transcritical bifurcation theorem by Castillo-Chavez and Song \cite{Castillo_Chavez} is used.

\begin{theorem} \label{Cas_th}
Consider a system of ordinary differential equations with a parameter $\phi$:
\begin{equation}\label{ODE}
\dot x(t)=f(x(t),\phi), \qquad f:\mathbb{R}^{n}\times\mathbb{R}\rightarrow \mathbb{R}^{n}  \quad \textrm{ and }  \quad f \in C^{2}(\mathbb{R}^{n}\times\mathbb{R}).
\end{equation}
It is assumed  that $0\in \mathbb{R}^{n}$ is an equilibrium for system (\ref{ODE}) for all values of the parameter $\phi$, that is
$f(0,\phi)= 0$ for all $\phi \in \mathbb{R}$.
Assume
\begin{enumerate}
\item $J=D_{x}f(0,0)=\Big( \frac{\partial f_{i}}{\partial x_{j}} (0,0)\Big)$ is the linearization of system (\ref{ODE}) around the equilibrium $0$ with $\phi$ evaluated at $0$. Zero is a simple eigenvalue of $J$ and all other eigenvalues of $J$ have negative real parts;
\item Matrix $J$ has a nonnegative right eigenvector $w$ and a left eigenvector $v$ corresponding to the zero eigenvalue.
\end{enumerate}

Let $f_k:\mathbb{R}^{n}\times\mathbb{R}\rightarrow \mathbb{R}$ be the $k$-th component of $f$ and
\begin{equation}\label{ab}
b=\displaystyle \sum _{k,i,j=1}^{n} v_k w_i w_j
\frac{\partial^{2}f_k}{\partial x_i \partial x_j}(0,0), \qquad
\quad d=\displaystyle \sum _{k,i=1}^{n} v_k w_i
\frac{\partial^{2}f_k}{\partial x_i \partial \phi}(0,0).
\end{equation}

The local dynamics of (\ref{ODE}) around $0$ are totally determined by $b$ and $d$ as follows. If $b<0$, $d>0$, then when $\phi$ changes from negative to
positive, with $|\phi|\ll 1$, $0$ changes its stability from stable to unstable. Correspondingly, a negative unstable equilibrium becomes positive
and locally asymptotically stable.
\end{theorem}

We note that the other cases concerning the signs of $b$ and $d$ are also considered in \cite{Castillo_Chavez}, however, here we need only this special case.

This theorem is applied for the following version of the compact pairwise model (\ref{alap_CP1})-(\ref{alap_CP5}), in which equations (\ref{alap_CP1}) and (\ref{alap_CP4}) are omitted and the parameter $\phi=\tau-\tau_{c}$ is introduced.
\begin{eqnarray}
 \dot{[I_l]}&=&-\gamma [I_l]+(\phi+\tau_c) n_l (N_l-[I_l]) \frac{[SI]}{nN-I_s}=:f_l, ~ l=1,\ldots, L \label{CP12''}\\
 \dot{[SI]}&=&\gamma ([II]-[SI])+(\phi+\tau_c) (nN-3[SI]-[II])[SI] \tilde{Q}-(\phi+\tau_c) [SI]=:f_{L+1}, \label{CP3''}\\
 \dot{[II]}&=&2(\phi+\tau_c) [SI]-2\gamma [II]+2(\phi+\tau_c) [SI]^2 \tilde{Q}=:f_{L+2} \label{CP5''},
\end{eqnarray}
where \[I_s=\sum_{l=1}^L n_l[I_l],\qquad \tilde{D}=\sum_{l=1}^L
n^{2}_l[I_l],\qquad \tilde{Q}=\frac{1}{(nN-I_s)^{2}} \big(\langle
n^2 \rangle N-nN+I_s-\tilde{D}\big).\]

\begin{theorem}
If $\tau_{c}<\tau$ with $|\tau-\tau_{c}|\ll 1$ then the endemic steady state is locally asymptotically stable in the compact pairwise model.
\end{theorem}

\noindent {\sc Proof.}

In system (\ref{CP12''})-(\ref{CP5''}) the disease free equilibrium point is $0=(0,\ldots,0,0,0) \in \mathbb{R}^{L+2}$, which is a steady state for all values of the parameter $\phi$. That is, $f(0,\phi)\equiv 0$ for $\phi \in (-\tau_c, +\infty)$, where $f:\mathbb{R}^{L+2}\times\mathbb{R}\rightarrow\mathbb{R}^{L+2}$
and $f:=(f_1, \ldots, f_L, f_{L+1},f_{L+2})$.

Similarly to the proof of Proposition \ref{prop2}, the Jacobian matrix of system (\ref{CP12''})-(\ref{CP5''}) at the equilibrium $0$ and $\phi=0$ is
$$
J_c=\begin{pmatrix}-\gamma I&R_c\\
                      0 & P_c
\end{pmatrix},
$$
where $0$ is a $2$-by-$L$ zero matrix and $I$ is a $L$-by-$L$ identity matrix. Moreover, $R_c$ denotes a $L$-by-$2$ matrix, in which the $l$-th entry in the first column is $\frac{\tau_c n_l N_l}{nN}$, $l=1,\ldots,L$, and all the elements in the second column are zeros. The matrix $P_c$ takes the form
$$
P_c=\begin{pmatrix}-\tau_c &\gamma\\
                      2\tau_c & -2\gamma
\end{pmatrix} \ .
$$
Using the block structure of $J_c$, it is obvious, that the multiplicity of the eigenvalue $-\gamma $ is $L$ and the other two eigenvalues are $-(2\gamma+\tau_c)$ and $0$. Hence, the first condition of Theorem \ref{Cas_th} is satisfied. The right eigenvector corresponding to the eigenvalue $0$ is
$$
w=(\ldots, \frac{\tau_{c}n_l N_l}{nN}, \ldots  ,\gamma, \tau_{c})^{T} \in \mathbb{R}^{(L+2)\times 1},
$$
and a left eigenvector is
$$
v=(0,\ldots, 0, 2, 1)\in \mathbb{R}^{1\times (L+2)},
$$
the coordinates of which are nonnegative, thus the second condition of Theorem \ref{Cas_th} is also fulfilled.

Now, let us calculate the values of $b$ and $d$. This can be achieved by calculating the second order partial derivatives of $f$ evaluated at the disease free equilibrium point with $\phi=0$. Simple calculations shows that
$$
\frac{\partial^{2} f_{L+1}}{\partial [I_l]\partial [SI]}(0,0)=-\frac{\tau_c}{nN}\frac{(n^{2}_l-n_l)n-(\langle n^2 \rangle -n)2n_l}{n}, \qquad
\frac{\partial^{2} f_{L+1}}{\partial [SI]^{2}}(0,0)=-\frac{6\tau_c}{nN}\frac{\langle n^2 \rangle-n}{n},
$$
$$
\frac{\partial^{2} f_{L+1}}{\partial [SI]\partial [II]}(0,0)=-\frac{\tau_c}{nN}\frac{\langle n^2 \rangle-n}{n}, \qquad
\frac{\partial^{2} f_{L+1}}{\partial [SI]\partial \phi}(0,0)=\frac{\langle n^2 \rangle-n}{n}-1,
$$
$$
\frac{\partial^{2} f_{L+2}}{\partial [SI]\partial [SI]}(0,0)=\frac{4\tau_c}{nN}\frac{\langle n^2 \rangle -n}{n}, \qquad
\frac{\partial^{2} f_{L+2}}{\partial [SI]\partial \phi}(0,0)=2.
$$
It is easy to see that the rest of the second derivatives are all zero. Thus, applying (\ref{ab}) we get
$$
b=2 \Big( \displaystyle \sum_{l=1}^{L} 2 \frac{\tau_c n_l N_l}{nN} \gamma \frac{-\tau_c}{nN}\frac{(n^{2}_l-n_l)n-(\langle n^2 \rangle -n)2n_l}{n} \Big) +
$$
$$
2\gamma^{2} \frac{-6\tau_c}{nN} \frac{\langle n^2 \rangle -n}{n} +4\gamma \tau_c \frac{-\tau_c}{nN} \frac{\langle n^2 \rangle -n}{n} + \gamma^2 \frac{4\tau_c}{nN} \frac{\langle n^2 \rangle -n}{n} ,
$$
$$
d=2 \gamma \left(\frac{\langle n^2 \rangle -n}{n} -1\right) + 2\gamma .
$$
It is obvious that $d= 2\gamma \frac{\langle n^2 \rangle -n}{n} >0.$ Simple algebra shows that the sum in $b$ can be rearranged as
$$
2 \displaystyle \sum_{l=1}^{L} 2 \frac{\tau_c n_l N_l}{nN} \gamma \frac{-\tau_c}{nN}\frac{(n^{2}_l-n_l)n-(\langle n^2 \rangle -n)2n_l}{n}=
$$
$$
-4 \frac{ \tau_c^2 \gamma}{n^2N^2} \frac{N}{n} \displaystyle \sum_{l=1}^{L}   n^{2}_l \frac{N_l}{N} \left(n_ln+n-2\langle n^2 \rangle \right)=
-4 \frac{ \tau_c^2 \gamma}{n^2N^2} \frac{N}{n}  \left[  n \langle n^3 \rangle + n  \langle n^2 \rangle -2\langle n^2 \rangle ^{2} \right],
$$
where $\langle n^3 \rangle=\displaystyle \sum_{l=1}^{L}  n^{3}_l \frac{N_l}{N}$ is a third moment of the degree distribution.

Using that $\tau_c (\langle n^2 \rangle - \langle n \rangle )=\gamma  \langle n \rangle $ and $n=\langle n \rangle$, the expression for $b$ can be easily transformed as follows.
\begin{align*}
b&=-4 \frac{ \tau_c^2 \gamma}{n^2N^2} \frac{N}{n} \big[  n \langle n^3 \rangle + n \langle n^2 \rangle -2\langle n^2 \rangle ^{2} \big ]  -8 \gamma^{2} \frac{\tau_c}{nN} \frac{\langle n^2 \rangle -n}{n}-4 \gamma  \frac{\tau_c^{2}}{nN} \frac{\langle n^2 \rangle -n}{n} \\
&=\gamma \tau_c \Bigg[-4 \frac{\gamma}{\langle n^2 \rangle -n} \frac{ 1 }{n^2N} \big[  n \langle n^3 \rangle + n \langle n^2 \rangle -2\langle n^2 \rangle ^{2} \big ]  -8 \frac{\gamma }{nN} \frac{\langle n^2 \rangle
-n}{n}-4 \frac{\gamma}{nN} \Bigg] \\
&=\frac{4\gamma^{2} \tau_c}{n^2N} \Big[- \frac{1}{\langle n^2 \rangle -n}\big[  n \langle n^3 \rangle + n \langle n^2 \rangle -2\langle n^2 \rangle ^{2} \big ] -2( \langle n^2 \rangle - n)-n \Big] \\
&=\frac{4\gamma^{2} \tau_c}{n^2N} \frac{1}{\langle n^2 \rangle -n}\left[- n \langle n^3 \rangle - n \langle n^2 \rangle +2\langle n^2 \rangle ^{2} -(2\langle n^2 \rangle- n)(\langle n^2 \rangle -n) \right]\\
&=\frac{4\gamma^{2} \tau_c}{n^2N} \frac{1}{\langle n^2 \rangle -n}\left[- n \langle n^3 \rangle  +2n \langle n^2 \rangle  -n^2\right]=\frac{4\gamma^{2} \tau_c}{nN} \frac{1}{\langle n^2 \rangle -n}\left[-   \langle n^3 \rangle  +2 \langle n^2 \rangle  -n\right]<0,
\end{align*}
where the last inequality follows from the fact that $2 \langle
n^2 \rangle< \langle n^3 \rangle +n$ holds for any degree
sequence, since $2n_l<n_l^2+1$ is true for $n_l>1$. $\Box$

\begin{remark}
We note that the endemic steady state exists mathematically also for $\tau<\tau_{c}$, however it has negative coordinates, hence it is not relevant from the epidemiological point of view. This is why it is claimed that the endemic steady state appears at the critical value of $\tau$. The general transcritical bifurcation theorem (Theorem \ref {Cas_th}) implies that it is unstable when $\tau<\tau_{c}$.
\end{remark}

\vspace*{0.5cm} \setcounter{equation}{0}
\section{Global stability of the disease-free steady state}

In this section, the global stability of the disease-free steady state is proved for
$$
\tau < \tau_c= \frac{\gamma n}{\langle n^{2}\rangle-n} := a\gamma
$$
when one of the following condition holds for the network.

\begin{itemize}
\item[(A1)] $(2+\sqrt{2})n\leq \langle n^{2}\rangle$

\item[(A2)] The network is bimodal, that is $L=2$.
\end{itemize}
In order to prove global stability, we introduce $\theta=\frac{[SI]}{S_s}$, that enables us to reduce the system to $L+1$ equations. Differentiating $\theta$ and using the differential equations (\ref{alap_CP1})-(\ref{alap_CP5}), we obtain the following system.

\begin{align}
\dot{[S_l]}(t)&=\gamma (N_l-[S_l](t))-\tau n_l [S_l](t) \theta(t), \quad l=1,\ldots, L, \label{CPtS}\\
\dot{\theta}(t)=&\gamma \frac{nN}{S_s(t)} (1-\theta(t)) -\gamma(1+\theta(t)) +\tau \left(\frac{D(t)}{S_s(t)}-2\right)\theta(t)(1-\theta(t)), \label{CPtheta}
\end{align}
where
\[S_s(t)=\sum_{l=1}^L n_l[S_l](t)=[SS](t)+[SI](t) \quad \mbox{ and } \quad D(t)=\sum_{l=1}^L n^{2}_l[S_l](t).\]
In order to prove the global stability result we will need two auxiliary Lemmas presented in the next subsection.

\subsection{Auxiliary results}

The following elementary statement will play a key role in the proof. Its proof is a simple exercise, it is recalled here only for completeness.

\begin{lemma}\label{lemma0}
Let $F:[0,1]\rightarrow [0,+\infty)$ be a continuous function, for which $F(x)<x$ for all $x\in(0,1]$ and $F(0)=0$.  Taking any initial point
$x_0\in[0,1]$, the sequence defined by $x_{n+1}:=F(x_n)$ is convergent and $\lim_{n\rightarrow + \infty}x_n=0$.
\end{lemma}

\noindent {\sc Proof.} The convergence of $(x_n)$ follows from the facts that it is decreasing and bounded below. Denoting its limit by $x^{*}$, the continuity of $F$ yields that $F(x^{*})=F(\lim_{n\rightarrow + \infty}x_n)=\lim_{n\rightarrow + \infty}F(x_n)=\lim_{n\rightarrow
+ \infty}x_{n+1}=x^{*}$ implying $x^{*}=0$.

$\Box$

We will use the following comparison theorem \cite{Hale}.

\begin{lemma} \label{comp_th}
Let $f(t,x)$ be a continuous function in $x$. Assume that
\begin{itemize}
\item the initial value problem $\dot{x}_2(t)=f(t,x_2(t))$, $x_2(0)=x_0$ has a unique solution for $t\in[0,T]$ and
\item $\dot{x}_1(t)\leq f(t,x_1(t))$ for $t\in[0,T]$ and $x_1(0)\leq x_0$.
\end{itemize}

Then $x_1(t)\leq x_2(t)$, for $t\in[0,T]$.
\end{lemma}

The global stability will be proved by the so-called monotone iteration technique \cite{Xinchu}, one step of which will be shown first.

\subsection{One step of the iteration}\label{it_lepes}

An iteration step consists of two parts. In the first one, it is shown that an upper bound of $\theta$ yields a lower bound of $S_l$. While in the second, it is proved that a lower bound of $S_l$ yields an upper bound of $\theta$. These statements will be shown in this subsection.

\begin{lemma} \label{lemma1}
Let $\tau < \tau_c=\gamma a$  and $x\in(0,1]$, $t_0\geq0$ be two numbers, such that $\theta(t)\in (0, x]$ for all $t\geq t_0$. Then there exists $t^*_1>t_0$, such that
\begin{equation}\label{Slbecsles}
\frac{N_l}{1+a n_l x} < [S_l](t), \quad \mbox{for } t>t^*_1 \quad \mbox{ and }  l=1,\ldots,L .
\end{equation}
\end{lemma}

\noindent {\sc Proof.}
Let us consider the differential equation (\ref{CPtS}) for a fixed $l=1,\ldots,L$:
$$
\dot{[S_l]}=\gamma N_l-[S_l] (\gamma+\tau n_l \theta),
$$
and the following constant-coefficient linear differential equation:
\begin{equation}\label{S_l_y}
\dot{y}=\gamma N_l-y (\gamma+\tau n_l x).
\end{equation}

Suppose that $[S_l](t_0)=y(t_0)$, then Lemma \ref{comp_th} implies $y(t)\leq [S_l](t)$ for all $t>t_0$, since $\gamma N_l-y (\gamma+\tau n_l x)\leq\gamma
N_l-[S_l] (\gamma+\tau n_l \theta)$.

The solution of the equation (\ref{S_l_y}) can be expressed as
$$
y(t)=\frac{\gamma}{\gamma+\tau n_l x}N_l+C e^{-(\gamma+\tau n_l x)t}, \textrm{ where } C= \big(y(t_0)-\frac{\gamma}{\gamma+\tau n_l x}N_l \big)e^{-(\gamma+\tau n_l x)t_0}.
$$
This leads to $\lim_{t\rightarrow + \infty}y(t)=\frac{\gamma}{\gamma+\tau n_l x}N_l$, consequently for all $\varepsilon>0$ there exists $t^*_1>t_0$,
such that if $t>t^*_1$, then we have $\frac{\gamma}{\gamma+\tau n_l x}N_l-\varepsilon<y(t)\leq [S_l](t)$.  Using $\tau <\gamma a$, one can choose $\varepsilon$ in such a way that $\varepsilon <\frac{N_l}{1+\frac{\tau}{\gamma} n_l x}-\frac{N_l}{1+a n_l x}$, leading to $\frac{N_l}{1+a n_l x} < [S_l](t)$ for $t>t^*_1$.

$\Box$

Let us turn to the second part of the iteration step, in which the previously obtained lower bound of $S_l$ yields an upper bound of $\theta$. First, we derive upper bounds for the coefficients in equation \eqref{CPtheta}, then we derive the upper bound for $\theta$.

\begin{lemma} \label{lemmaSs}
Let $\tau < \tau_c=\gamma a$  and assume that there is an $x\in(0,1]$ and $t^*_1>0$, such that the lower bound in \eqref{Slbecsles} holds. Then introducing $B=\frac{\langle n^{2}\rangle}{\langle n^{2}\rangle -n}$, we have
\begin{equation}\label{NnperSsbecsles}
\frac{nN}{S_s(t)} \leq 1+Bx, \quad \mbox{for } t>t^*_1.
\end{equation}
\end{lemma}

\noindent {\sc Proof.}
The lower bound in \eqref{Slbecsles} yields
$$
g(x):=\sum_{l=1}^L \frac{n_l N_l}{1+a n_l x}<S_s.
$$
We will show that $\frac{1}{1+Bx}\leq \frac{g(x)}{nN}$ by using Jensen's inequality. Let us consider the convex function $f:[0,+\infty)\rightarrow \mathbb{R}$,
$f(x)=\frac{1}{1+x}$ and let $c_l=\frac{n_l N_l}{nN}$, $y_l=a n_l x$, $l=1,\ldots,L$ for which $\sum_{l=1}^L c_l=\sum_{l=1}^L \frac{n_l N_l}{nN}=1$ holds. Hence, applying Jensen's inequality, we get
$$
\sum_{l=1}^L \frac{n_l N_l}{nN}\frac{1}{1+a n_l x}=\sum_{l=1}^L c_l f(y_l) \geq f(\sum_{l=1}^L c_l y_l)= \frac{1}{1+\sum_{l=1}^L \frac{n_l N_l}{nN}a n_l x}=\frac{1}{1+\frac{ax}{n}\sum_{l=1}^L \frac{n_l N_l}{N} n_l}=
$$
$$
=\frac{1}{1+\frac{ax}{n}\langle n^{2}\rangle}=\frac{1}{1+\frac{n}{\langle n^{2}\rangle-n}\frac{x}{n}\langle n^{2}\rangle}=\frac{1}{1+\frac{\langle n^{2}\rangle}{\langle n^{2}\rangle -n}x}=\frac{1}{1+Bx}.
$$
$\Box$

\begin{lemma} \label{lemmaDperSs}
Let $\tau < \tau_c=\gamma a$  and assume that there is an $x\in(0,1]$ and $t^*_1>0$, such that the lower bound in \eqref{Slbecsles} holds. Then we have
\begin{equation}\label{DperSsbecsles}
\frac{D(t)}{S_s(t)} \leq \frac{\langle n^{2}\rangle}{n}(1+Bx), \quad \mbox{for } t>t^*_1.
\end{equation}
Moreover, if the matrix is bimodal, that is $L=2$, then we have the following alternative upper bound
\begin{equation}\label{DperSsbecslesbimod}
\frac{D(t)}{S_s(t)} \leq \frac{ n^{2}_1N_1+ (1+a n_1 x)n^{2}_2 N_2}{ n_1N_1+ (1+a n_1 x)n_2 N_2}, \quad \mbox{for } t>t^*_1.
\end{equation}
\end{lemma}

\noindent {\sc Proof.}
First, (\ref{DperSsbecsles}) is proved by using $\sum\limits_{l=1}^L n_l^2 [S_l] \leq \sum\limits_{l=1}^L n_l^2 N_l = \langle n^{2}\rangle N$ and Lemma \ref{lemmaSs} as
$$
\frac{D(t)}{S_s(t)}\leq \frac{\langle n^{2}\rangle N}{S_s(t)}\leq
\frac{\langle n^{2}\rangle}{n}(1+Bx).
$$

Secondly, we verify (\ref{DperSsbecslesbimod}) by calculating the maximum of the function $f(x_1, x_2)=\frac{ n^{2}_1 x_1+ n^{2}_2 x_2}{
n_1 x_1+ n_2 x_2}$  on the rectangle $T:=[\frac{1}{1+a n_1 x}N_1, N_1]\times [\frac{1}{1+a n_2 x}N_2,N_2]$. It is $\underline{x}=(\frac{1}{1+a n_1 x}N_1,N_2)$, assuming $n_1<n_2$, hence
$$\frac{D(t)}{S_s(t)}=\frac{ n^{2}_1[S_1](t)+n^{2}_2[S_2](t)}{
n_1[S_1](t)+n_2[S_2](t)}\leq \frac{ n^{2}_1 \frac{1}{1+a n_1 x}N_1+ n^{2}_2 N_2}{ n_1
\frac{1}{1+a n_1 x}N_1+ n_2 N_2}= \frac{ n^{2}_1N_1+ (1+a n_1
x)n^{2}_2 N_2}{ n_1N_1+ (1+a n_1 x)n_2 N_2}.$$

$\Box$

Now we have the upper estimates for the coefficients in equation \eqref{CPtheta}, hence we are ready to derive the upper bound for $\theta$.

\begin{lemma} \label{lemmatheta}
Let $\tau < \tau_c=\gamma a$  and assume that there is an $x\in(0,1]$ and $t^*_1>0$, such that the lower bound in \eqref{Slbecsles} holds. Moreover, assume that $\frac{D(t)}{S_s(t)} \leq b(x)$ with one of the functions obtained in Lemma \ref{lemmaDperSs}. Let us introduce the quadratic polynomial
$$
p_x(z)=\gamma (1+Bx) (1-z) -\gamma(1+z) +\gamma a \left(b(x)-2\right)z(1-z), \quad z \in [0,1].
$$
Then $p_x(0)>0$, $p_x(1)<0$ and $p_x$ has a unique root $z^*(x)\in (0,1)$. Moreover, for any $\overline{z}\in (z^*(x),1)$ there is a number $t_1>t^*_1$ such that the solution of \eqref{CPtheta} satisfies $\theta(t)< \overline{z}$ for $t>t_1$.
\end{lemma}

\noindent {\sc Proof.}
It is easy to see, that $p_x(0)>0$ and $p_x(1)<0$ hold, since $p_x(0)=\gamma B x >0$ and $p_x(1)=-2\gamma<0$, yielding that $p_x$ has a unique root $z^*(x)\in (0,1)$. In order to prove the second part of the statement, let us consider the following autonomous differential equation:
\begin{equation}
\dot{y}(t)=p_x(y(t))   \textrm{ with the initial condition } 0<y(t^*_1)=y_0\leq x. \label{de_px}
\end{equation}

It is clear, that the solution of (\ref{de_px}), denoted by $y(t)$, converges to $z^{*}(x)$ as $t \rightarrow +\infty$, i.e.
$\forall \varepsilon
>0$ there exists $t_1>t^*_1$, such that, if $t>t_1$, then $y(t)<z^{*}(x)+\varepsilon$ holds.

Now, let us denote the right hand side of equation (\ref{CPtheta}) by
$$q_t(\theta):=\gamma \frac{nN}{S_s(t)} (1-\theta) -\gamma(1+\theta) +\tau \left(\frac{D(t)}{S_s(t)}-2\right)\theta(1-\theta).$$
Thus $\theta$ is the solution of the initial value problem
\begin{equation}
\dot{\theta}(t)=q_t(\theta(t))   \textrm{ and } \theta(t^*_1)=y_0. \label{fam_de}
\end{equation}
According to Lemmas \ref{lemmaSs} and \ref{lemmaDperSs} we have $q_t(z)\leq p_x(z)$, $\forall z \in [0,1]$. Hence applying the comparison result, Lemma \ref{comp_th} to the initial value problems (\ref{de_px}) and (\ref{fam_de}) we get that $\theta(t)\leq y(t)$, $\forall t >t^*_1 $, yielding that $\theta(t) <z^{*}(x)+\varepsilon$ holds, if $t>t_1$. By choosing $\varepsilon:= \overline{z}-z^*(x)$, the result is obtained.

$\Box$

In order to get the global stability by using the monotone iteration technique, we need to prove that $\theta$ gets closer to zero in each iteration step, that can be proved by showing that $z^*(x)<x$.

\begin{lemma} \label{lemmazstar}
Let $\tau < \tau_c=\gamma a$  and assume that there is an $x\in(0,1]$ and $t^*_1>0$, such that the lower bound in \eqref{Slbecsles} holds. Let us choose the function $b$ in the polynomial $p_x$ as follows.
\begin{enumerate}
\item[(i)] If $(2+\sqrt{2})n\leq \langle n^{2}\rangle$, then let $b(x)=\frac{\langle n^{2}\rangle}{n}(1+Bx)$.

\item[(ii)] If he network is bimodal, that is $L=2$, then let $b(x)=\frac{ n^{2}_1N_1+ (1+a n_1 x)n^{2}_2 N_2}{ n_1N_1+ (1+a n_1 x)n_2 N_2}$.
\end{enumerate}
Then for both cases we have that the root $z^*(x)$ of $p_x$ satisfies $z^*(x)<x$ for all $x\in (0,1]$.
\end{lemma}

\noindent {\sc Proof.}
In both cases, it is enough to prove that $p_x(x)<0$ for all $x \in (0,1]$, because this implies that for the root of $p_x$ the inequality $z^*(x)<x$ holds.

In case (i), in order to confirm
$$
p_x(x)=\gamma (1+Bx) (1-x) -\gamma(1+x) +\gamma a \left(\frac{\langle n^{2}\rangle}{n}(1+Bx)-2\right)x(1-x)<0,
$$
it is enough to show that the coefficients of $p_x$, i.e. the derivatives of $p_x$ at $0$ are not positive, and the leading coefficient is negative, since $p_x$ is a cubic polynomial. It is obvious, that $p_x(0)=0$ and $p'_x(0)=0$ can be easily seen. For the second derivative we have that
$$
p''_x(0)=-\frac{2\gamma}{(\langle n^{2}\rangle -n)^{2}}(\langle n^{2}\rangle^{2}-4\langle n^{2}\rangle n+2 n^{2})\leq 0 \Leftrightarrow (2+\sqrt{2})n\leq \langle
n^{2}\rangle,
$$
furthermore
$$
p^{(3)}_x(0)=-\frac{6\gamma\langle n^{2}\rangle}{(\langle n^{2}\rangle -n)^{2}}<0.
$$

In case (ii), we need to verify
$$
p_x(x)=\gamma(1+Bx)(1-x)-\gamma(1+x)+ \gamma a \Big( \frac{ n^{2}_1N_1+ (1+a n_1 x)n^{2}_2 N_2}{ n_1N_1+ (1+a n_1 x)n_2 N_2}-2\Big)x(1-x)<0.
$$
Multiplying $p_x(x)$ with the positive denominator $U(x):= n_1N_1+ (1+a n_1 x)n_2 N_2$, we get
$$
r(x):=\gamma(1+Bx)(1-x)U(x)-\gamma(1+x)U(x)+ \gamma a\Big( n^{2}_1N_1+ (1+a n_1 x)n^{2}_2 N_2-2U(x)\Big)x(1-x).
$$
Now, it is enough to see that $r(x)<0$, if $x \in (0,1]$. Since $r$ is a cubic polynomial, the proof will be similar to the one above. It is obvious, that $r(0)=0$, and $r'(0)=0$ hold. The second derivative of $r$ at $0$ is:
$$
r''(0)=-\frac{2\gamma n}{N(\langle n^{2}\rangle -n)^{2}}\big(2N^{2}(\langle n^{2}\rangle -n)^{2}+N(\langle n^{2}\rangle -n n_2)N_2 n_1 n_2 \big).
$$
We will show that the expression $V:=2N^{2}(\langle n^{2}\rangle -n)^{2}+N(\langle n^{2}\rangle -n n_2)N_2 n_1 n_2$
is nonnegative. For this, we use $N\langle n^{2}\rangle=n^{2}_1 N_1+n^{2}_2 N_2$ and $Nn=n_1 N_1+n_2 N_2$ to yield
$$
N^{2}(\langle n^{2}\rangle -n)^{2}=(n^{2}_1 N_1+n^{2}_2 N_2-n_1 N_1-n_2 N_2)^{2},
$$
$$
N(\langle n^{2}\rangle -n n_2)=(n_1-n_2)n_1N_1.
$$
A
simple calculation shows that
$$
V=2n^2_1N^2_1(n_1-1)^2+2n^2_2N^2_2(n_2-1)^2+N_1N_2n_1n_2(n_1-2)^2+N_1N_2n_1n^2_2(3n_1-4).
$$
It is easy to see, that $V\geq 0$, if $n_2>n_1\geq 2$. In the case when $n_1=1$, we can see that:
$$
V=2n^2_2N^2_2(n_2-1)^2+N_1N_2n_2+N_1N_2n^2_2(-1)=N_2n_2(n_2-1)(2n_2N_2(n_2-1)-N_1)\geq
$$
$$
\geq N_2n_2(n_2-1)(2n_2N_2-N_1)>0,
$$
since in case of a connected network a node of degree $1$ should join to a node of degree $n_2$, consequently the total number of stubs starting from nodes of degree $n_2$ is not less than the total number of stubs starting from nodes of degree $1$, namely $ n_2 N_2\geq N_1$.
Finally, the third derivative of $r$ at $0$ is:
$$
r^{(3)}(0)=-\frac{6\gamma n_2N_2n n_1}{(\langle n^{2}\rangle-n)^2}(\langle n^{2}\rangle+n n_2-2n)<0,
$$
taking into account the inequalities $\langle n^{2}\rangle > n$ and $n n_2>n$, when $n_2>n_1$.

$\Box$

Now, we are ready to prove the main result.

\subsection{Proof of the global stability of the disease free equilibrium}

In this subsection we prove the following main theorem.

\begin{theorem}
Let $\tau < \tau_c=\gamma a$  and assume that the network satisfies either assumption (A1) or (A2). Then the disease-free equilibrium is globally asymptotically stable. That is the solutions of (\ref{CPtS})-(\ref{CPtheta}) starting from any initial condition converge to the disease free steady state: $\lim_{t\rightarrow + \infty}S_l(t)=N_l$, $l=1,\ldots,L$ and $\lim_{t\rightarrow + \infty}\theta(t)=0$.
\end{theorem}

\noindent {\sc Proof.}

We apply the monotone iteration technique, the idea of which is to define a decreasing sequence $(x_n)$ tending to zero and then show that there is an increasing sequence $t_n$, such that $\theta(t)<x_{n}$ when $t>t_n$. This proves that $\lim_{t\rightarrow + \infty}\theta(t)=0$, which implies $\lim_{t\rightarrow + \infty}S_l(t)=N_l$ by using Lemma \ref{lemma1}.

Let us define the sequence $(x_n)$ by $x_0=1$ and $x_{n+1}=F(x_n)$, where $F:[0,1]\rightarrow [0,+\infty)$ is defined as
$$
F(x):=\frac{x+z^{*}(x)}{2},
$$
where $z^{*}(x)$ is the unique root of the polynomial $p_x$ in the interval $[0,1]$. Besides that, let us extend the function $F$ continuously to the closed interval $[0,1]$ by defining $F(0):=\lim_{x\rightarrow 0}\frac{x+z^{*}(x)}{2}=0$.

According to Lemma \ref{lemmazstar} we have $z^*(x)<x$, hence $F(x)<x$ holds for all $x\in (0,1]$. Thus Lemma \ref{lemma0} implies that $(x_n)$ tends to zero.

Now let us create the iteration. The initiation of the iteration is $x_0=1$ and $t_0=0$. Then $\theta(t)\in (0,x_0]$ holds for all $t>t_0$. Applying Lemma \ref{lemma1}, we get that there exists $t_1^*>t_0$, such that \eqref{Slbecsles} holds when $t>t_1^*$. Then we can take $\overline{z}=x_1$ in Lemma \ref{lemmatheta} and obtain that there exists $t_1>t_1^*$ such that $\theta(t)\in (0,x_1]$ holds for all $t>t_1$. The next steps of the iteration are made in the same way, that is $t_n$ is a value, for which $\theta(t)\in (0,x_n]$ holds for all $t>t_n$. This completes the proof of the theorem.
$\Box$

\section{Discussion}

\begin{figure}[h!]
   \centerline{ \includegraphics[width=14cm,height=7cm]{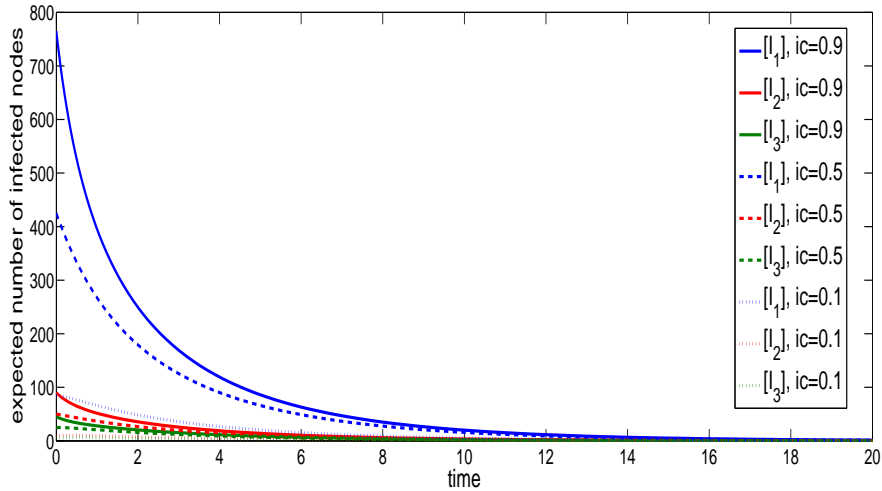}}
\caption{ Case of the globally stable disease-free equilibrium: Time evolution of the expected number of the infected nodes $[I_1]$, $[I_2]$, $[I_3]$ of degree $n_1=2$, $n_2=3$, $n_3=4$ respectively, started with $90$, $50$, $10$ randomly chosen infected nodes. The parameters are: $N=1000$, $N_1=850$, $N_2=100$, $N_3=50$, $\gamma=1$, $\tau=0.5$, $\tau_c=0.7586$. }
\label{cpwfig_0a}
\end{figure}

The global behaviour of the compact pairwise model of SIS epidemic propagation on a network, system (\ref{alap_CP1})-(\ref{alap_CP5}), was studied. We proved that transcritical bifurcation occurs at $\tau = \tau _c =  \frac{\gamma n}{\langle n^{2}\rangle-n}$. For subcritical values of $\tau$ the disease-free steady state is stable, while for supercritical values a unique stable endemic equilibrium appears. We also studied the global stability of the system. For subcritical values of $\tau$ we proved the global stability of the disease-free steady state under assumption (A1) and (A2). We note that these assumptions cover a wide class of networks. For example, it is easy to show that if each node has at least degree 4, then (A1) holds. However, there are graphs which satisfy neither (A1) nor (A2). An example is a network with the parameters: $N_1=850$, $N_2=100$, $N_3=50$, $n_1=2$, $n_2=3$, $n_3=4$, which is not bimodal and elementary calculation shows that (A1) is violated. Despite of this fact, the disease-free steady state is globally stable for subcritical $\tau$ values as Figure \ref{cpwfig_0a}
shows. We checked that the number of infected nodes tends to zero starting from  different initial conditions. Extensive numerical experiments show that the disease-free steady state is globally stable for any subcritical value of $\tau$, i.e. the assumptions (A1) and (A2) are not necessary.

\begin{figure}[h!]
   \centerline{ \includegraphics[width=14cm,height=7cm]{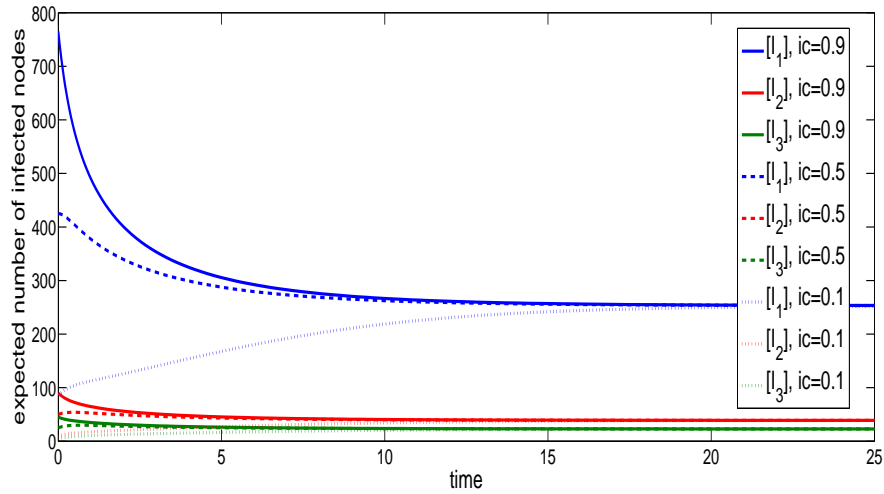}}
\caption{ Case of the globally stable endemic equilibrium: Time evolution of the expected number of the infected nodes $[I_1]$, $[I_2]$, $[I_3]$ of degree $n_1=2$, $n_2=3$, $n_3=4$ respectively, started with $90$, $50$, $10$ randomly chosen infected nodes. The parameters are: $N=1000$, $N_1=850$, $N_2=100$, $N_3=50$, $\gamma=1$, $\tau=1$, $\tau_c=0.7586$.}
\label{cpwfig_0b}
\end{figure}

We investigated the global stability of the endemic equilibrium for supercritical values of $\tau$ numerically and found that it is globally stable. An example is presented in Figure \ref{cpwfig_0b}, where the time dependence of the number of infected nodes is shown for different initial conditions when $\tau >\tau_c$. The analytic study of the global stability of the endemic steady state will be the subject of future work.

\section*{Acknowledgement}

Péter L. Simon acknowledges support from Hungarian Scientific Research Fund, OTKA,
(grant no. 115926).

The project has been supported by the European Union, co-financed by the European Social Fund (EFOP-3.6.3-VEKOP-16-2017-00002).

\end{document}